\documentclass{article}
\usepackage{graphicx}
\usepackage{amsmath}
\usepackage{amsfonts}
\usepackage{amssymb}

\newtheorem{theorem}{Theorem}

\newtheorem{corollary}{Corollary}

\newtheorem{definition}{Definition}

\newtheorem{lemma}{Lemma}

\newtheorem{proposition}{Proposition}

\newcommand{\K}{$\Bbb K$}
\newcommand{\U}{\mathcal{U}(\frak{g})}
\newcommand{\w}{\omega}

\begin{document}

\title{ Valued deformations of algebras}

\author{Michel Goze \thanks{%
corresponding author: e-mail: M.Goze@uha.fr} \and Elisabeth
Remm \thanks{E.Remm@uha.fr. Partially supported by a grant from the Institut of mathematics Simon Stoilow of the Romanian academy. Bucharest }\\
Universit\'{e} de Haute Alsace, F.S.T.\\
4, rue des Fr\`{e}res Lumi\`{e}re - 68093 MULHOUSE - France}
\date{}
\maketitle

\begin{abstract}
We develop the notion of  deformations  using a
valuation ring as ring of coefficients. This permits to consider in particular the classical Gerstenhaber
deformations of associative or Lie algebras as infinitesimal deformations and to solve the equation of deformations in a polynomial frame. We consider also the deformations of the enveloping algebra of
a rigid Lie algebra and we define valued deformations for some classes of non associative algebras.
\end{abstract}

\medskip

Table of contents :

1. Valued deformations of Lie algebras

2. Decomposition of valued deformations

3. Deformations of the enveloping algebra of a rigid Lie algebra

4. Deformations of non associative algebras

\bigskip

\section{Valued deformations of Lie algebras}

\subsection{Rings of valuation}
We recall briefly the classical notion of ring of valuation. Let $\Bbb{F}$ be a (commutative) field and $A$ a subring of $\Bbb{F}$.
We say that $A$ is a ring of valuation of $\Bbb{F}$ if $A$ is a local integral domain satisfying:
$$ \mbox{ If} \ x \in \Bbb{F} - A, \quad \mbox{then} \quad x^{-1} \in \frak{m}.$$
where $\frak{m}$ is the maximal ideal of $A$.

A ring $A$ is called ring of valuation if  it is a ring of valuation of its field of fractions.

\medskip

\noindent Examples : Let $\Bbb{K}$ be a commutative field of characteristic $0$. The ring of formal series $\Bbb{K}[[t]]$ is a valuation ring. 
On other hand the ring $\Bbb{K}[[t_1,t_2]]$ of two (or more) indeterminates
is not a valuation ring.

\subsection{Versal deformations of Fialowski [F] }
Let $\frak{g}$ be a $\mathbb{K}$-Lie algebra and $A$ an unitary commutative local $\mathbb{K}$-algebra. The tensor product 
$ \frak{g} \otimes A $ is naturally
endowed with a Lie algebra structure :
$$[ X \otimes a, Y \otimes b]=[X,Y] \otimes ab.$$
If $ \epsilon : A \longrightarrow \mathbb{K}$, is an unitary augmentation with kernel the maximal ideal $\frak{m}$, a deformation $\lambda $ of $\frak{g}$
 with base $A$ is a Lie algebra structure on $\frak{g} \otimes A$ with bracket $[,]_{\lambda}$ such that
$$ id \otimes \epsilon : \frak{g} \otimes A \longrightarrow
\frak{g} \otimes \Bbb{K}$$
is a Lie algebra homomorphism. In this case the bracket $[,]_{\lambda}$ satisfies
$$[X \otimes 1,Y \otimes 1]_{\lambda}=[X,Y] \otimes 1 + \sum Z_i \otimes a_i$$
where $a_i \in A$ and $X,Y,Z_i \in \frak{g}.$
Such a deformation is called infinitesimal if the maximal ideal $\frak{m}$ satisfies $\frak{m}^2 =0.$ An interesting example is described in [F].
If we consider the commutative algebra $A= \Bbb{K} \oplus
(H^2(\frak{g}, \frak{g}))^*$ (where $^*$ denotes the dual as vector space) such that $dim(H^2) \leq \infty$, the deformation
with base
$A$ is an infinitesimal deformation (which plays the role of an universal deformation).

\subsection{Valued deformations of Lie algebra}
Let $\frak{g}$ be a $\mathbb{K}$-Lie algebra and $A$ a commutative $\mathbb{K}$-algebra of valuation. Then $\frak{g} \otimes A$ is a $\mathbb{K}$-Lie algebra.
We can consider this Lie algebra as an $A$-Lie algebra. We denote this last by $\frak{g}_A$. If $dim_{\mathbb{K}}(\frak{g})$ is finite then
$$dim_{A}(\frak{g}_A)=dim_{\mathbb{K}}(\frak{g}).$$
As the valued ring $A$ is also a  $\mathbb{K}$-algebra  we have a natural embedding of the $\mathbb{K}$-vector space $\frak{g}$ into
the free
$A$-module $\frak{g}_A$. Without loss of generality we can consider this embedding to be the identity map.
\begin{definition}
Let $\frak{g}$ be  a $\mathbb{K}$-Lie algebra and $A$ a commutative $\mathbb{K}$-algebra of valuation such that the residual field
$\frac{A}{\frak{m}}$ is isomorphic to $\mathbb{K}$ (or to a subfield of $\mathbb{K}$). A valued deformation of $\frak{g}$ with base $A$ is a $A$-Lie algebra $\frak{g}'_A$
such that the underlying $A$-module of $\frak{g}'_A$ is $\frak{g}_A$ and that
$$[X,Y]_{\frak{g}'_A} -[X,Y]_{\frak{g}^{ \,} _A} $$ is in the $\frak{m}$-quasi-module $\frak{g} \otimes \frak{m}$ where $\frak{m}$ is the maximal ideal of
$A$.
\end{definition}
The classical notion of deformation studied by Gerstenhaber ([G]) is a valued deformation.
In this case $A=\Bbb{K}[[t]]$ and the residual field of $A$ is isomorphic to \K \, . Likewise a versal deformation is a valued deformation.
The algebra $A$ is in this case the finite dimensional \K-vector space $\Bbb{K} \oplus (H^2( \frak{g},\frak{g}))^*$ where $H^2$ denotes the
second Chevalley cohomology
group of $\frak{g}$. The algebra law  is given by
$$(\alpha _1, h_1).(\alpha _2, h_2)=(\alpha _1.\alpha _2, \alpha _1.h_2 + \alpha _2.h_1).$$
It is a local field with maximal ideal $\{0\} \oplus (H^2)^*$. It is also a valuation field because we can endowe this algebra with a field structure,
the inverse of $(\alpha, h)$ being $((\alpha)^{-1}, -(\alpha)^{-2}h)$.

\section{Decomposition of valued deformations}

In this section we show that every valued deformation can be decomposed in a finite sum (and not as a serie) with pairwise comparable infinitesimal coefficients (that is 
in $\frak{m}$). The interest of this decomposition is to avoid the classical problems of convergence.  
\subsection{Decomposition in $\frak{m} \times \frak{m}$}
Let $A$ be a valuation ring satisfying the conditions of definition
1. Let us denote by $\cal{F}_A$ the field of fractions of $A$ and
$\frak{m}^2$ the catesian product $\frak{m} \times \frak{m}$ .
Let $(a_1,a_2) \in \frak{m}^2$ with $a_i \neq 0$ for $i=1,2$.

\smallskip

\noindent i)
Suppose that $a_1.a_2^{-1} \in A$ and $a_2a_1^{-1} \in A$. Let be $\alpha = \pi (a_1.a_2^{-1})$ where $\pi$ is the canonical projection on $\frac{A}{\frak{m}}$.
Clearly, there exists a global section
$s: \Bbb{K} \rightarrow A$ which permits to identify $\alpha$ with $s(\alpha)$ in $A$. Then
$$a_1.a_2^{-1}=\alpha + a_3$$
with $a_3 \in \frak{m}$.
Then if $a_3 \neq 0$, 
$$(a_1,a_2)=(a_2(\alpha + a_3),a_2)=a_2(\alpha,1)+a_2a_3(0,1).$$
If $\alpha \neq 0$ we can also write
$$(a_1,a_2)=aV_1+abV_2$$
with $a,b \in \frak{m}$ and $V_1,V_2$ linearly independent in $\Bbb{K}^2$. If $\alpha = 0$ then $a_1.a_2^{-1} \in \frak{m}$ and $a_1=a_2a_3$.
We have
$$(a_1,a_2)=(a_2 a_3,a_2)=ab(1,0)+a(0,1).$$
So in this case, $V_1=(0,1)$ and $V_2=(1,0)$. If $a_3 =0$ then
$$a_1a_2^{-1}=\alpha$$
and 
$$(a_1,a_2)=a_2(\alpha,1)= aV_1.$$

This correspond to the previous decomposition but with $b=0$.

\smallskip

\noindent ii) If
$a_1.a_2^{-1} \in {\cal{F}}_A-A$, then $a_2.a_1^{-1} \in \frak{m}$. We put in this case  $a_2.a_1^{-1}=a_3$ and we have
$$(a_1,a_2)=(a_1,a_1.a_3)=a_1(1,a_3)=a_1(1,0)+a_1a_3(0,1)$$
with $a_3 \in \frak{m}$. Then, in this case the point $(a_1,a_2)$  admits the following decomposition :
$$(a_1,a_2)=aV_1 + ab V_2$$
with $a,b \in \frak{m}$ and $V_1,V_2$ linearly independent in $\Bbb{K}^2$. Note that this case corresponds to the previous but with $\alpha =0$.

\smallskip

\noindent iii)
If $a_1.a_2^{-1} \in A$ and $a_2.a_1^{-1}\notin A$
then as $a_2.a_1^{-1} \in \mathcal{F}_A-A$, $a_1 a_2^{-1} \in \frak{m}$ and we find again the precedent case with
$\alpha=0$.
Then we have proved
\begin{proposition}
For every point $(a_1,a_2) \in \frak{m}^2$, there exist lineary independent
vectors $V_1$ and $V_2$ in the $\mathbb{K}$-vector space $\mathbb{K} ^2$ such that
$$(a_1,a_2)=aV_1+abV_2$$
for some $a,b \in \frak{m}$. 
\end{proposition}
Such decomposition est called of length $2$ if $b \neq 0$. If not it is called of length $1$. 

\subsection{Decomposition in $\frak{m}^k$}
Suppose that $A$ is valuation ring satisfying the hypothesis of Definition 1. Arguing as before, we can conclude
\begin{theorem}
For every $(a_1,a_2,...,a_k) \in \frak{m}^k$ there exist $h \ (h \leq k$) independent vectors $V_1,V_2,..,V_h$
whose components are in $\mathbb{K}$ and elements
$b_1,b_2,..,b_h \in \frak{m}$ such that
$$(a_1,a_2,...,a_k)=b_1V_1+b_1b_2V_2+...+b_1b_2...b_hV_h.$$
\end{theorem}
The parameter $h$ which appears in this theorem is called the length of the decomposition. This parameter can be different to $k$. It corresponds
to the dimension of the smallest $\mathbb{K}$-vector space $V$ such that $(a_1,a_2,...,a_k) \in V \otimes \frak{m}$. 

\smallskip

\noindent If the coordinates $a_i$ of the vector $(a_1,a_2,...,a_k)$ are in $A$ and not necessarily in its maximal ideal, then writing $a_i=\alpha _i +a'_i$
with $\alpha _i \in \mathbb{K} $ and $a_i' \in \frak{m}$,
we decompose
$$(a_1,a_2,...,a_k)=(\alpha _1,\alpha _2,...,\alpha _k) +(a'_1,a'_2,...,a'_k)$$
and we can apply Theorem 1 to the vector $(a'_1,a'_2,...,a'_k)$.

\subsection{Uniqueness of the decomposition}
Let us begin by a technical lemma.
\begin{lemma}
Let $V$ and $W$ be two vectors with components in the valuation ring $A$. There exist $V_0$ and $W_0$ with components in $\mathbb{K}$
such that $V=V_0+ V'_0$ and $W=W_0 + W'_0$ and the components of $V'_0$ and $W'_0$ are in the maximal ideal $\frak{m}$. Moreover if
the vectors $V_0$ and $W_0$ are linearly independent then $V$ and $W$ are also independent.
\end{lemma}
{\it Proof.} The decomposition of the two vectors $V$ and $W$ is evident. It remains to prove that the independence of the vectors $V_0$
and $W_0$ implies those of $V$ and $W$. Let $V,W$ be two vectors with components in $A$ such that $\pi (V)=V_0$ and $\pi (W)=W_0$ are independent. Let us suppose that $$xV+yW=0$$ with $x,y \in A$. One of the coefficients $xy^{-1}$ or $yx^{-1}$ is not in $\frak{m}$. Let us suppose that $xy^{-1} \notin \frak{m}$. If $xy^{-1} \notin A$ then $x^{-1}y \in \frak{m}$. Then
$ xV+yW=0$ is equivalent to $V+x^{-1}yW=0$. This implies that $\pi (V)=0$ and this is impossible. Then $xy^{-1} \in A-\frak{m}$. Thus if there exists a linear relation between $V$ and $W$, there exists a linear relation with coefficients in $A-\frak{m}$. We can suppose that $xV+yW=0$ with $x,y \in A-\frak{m}$. As
$V=V_0+V'_0, \, W=W_0+W'_0$ we have
$$\pi (xV+yW)=\pi (x)V_0+\pi (y)W_0=0.$$
Thus $\pi(x)=\pi(y)=0$. This is impossible and the vectors $V$ and $W$ are independent as soon as $V_0$ and $W_0$ are independent vectors. $\quad \Box$

\bigskip

Let $(a_1,a_2,...,a_k)=b_1V_1+b_1b_2V_2+...+b_1b_2...b_hV_h$ and $(a_1,a_2,...,a_k)=c_1W_1+c_1c_2W_2+...+c_1c_2...c_sW_s$ be two decompositions
of the vector $(a_1,a_2,...,a_k)$. Let us compare the coefficients $b_1$ and $c_1$. By hypothesis $b_1c_1^{-1}$ is in $A$ or the inverse is in $\frak{m}$. Then we can suppose that $b_1c_1^{-1} \in A$.
As the residual field is a subfield of \K \, , there exists $\alpha \in \frac{A}{\frak{m}}$ and $c_1 \in \frak{m}$ such that
$$b_1c_1^{-1} = \alpha + b_{11}$$
thus $b_1=\alpha c_1+ b_{11}c_1$. Replacing this term in the decompositions we obtain
$$
\begin{array}{l}
(\alpha c_1+ b_{11}c_1)V_1+(\alpha c_1+ b_{11}c_1)b_2V_2+...+(\alpha c_1+ b_{11}c_1)b_2...b_hV_h \\ =c_1W_1+c_1c_2W_2+...+c_1c_2...c_sW_s.
\end{array}
$$
Simplifying by $c_1$, this expression is written
$$\alpha V_1+ m_1 = W_1 + m_2$$
where $m_1,m_2$ are vectors with coefficients $\in \frak{m}$. From Lemma 1, if $V_1$ and $W_1$ are linearly independent, as its coefficients are in the residual field,
the vectors $\alpha V_1+ m_1$ and $W_1 + m_2$ would be also linearly independent ($\alpha \neq 0$). Thus $W_1=\alpha V_1$. One deduces
$$b_1V_1+b_1b_2V_2+...+b_1b_2...b_hV_h=c_1(\alpha V_1)+c_1b_{11}V_1+c_1b_{12}V_2+...+c_1b_{12}b_3...b_hV_h,$$ with $b_{12}=b_2(\alpha+b_{11}).$
Then
$$b_{11}V_1+b_{11}b_{12}V_2+...+b_{11}b_{12}b_3...b_hV_h=c_2W_2+...+c_2...c_sW_s.$$
Continuing this process by induction we deduce the following result
\begin{theorem}
Let be $b_1V_1+b_1b_2V_2+...+b_1b_2...b_hV_h$ and $c_1W_1+c_1c_2W_2+...+c_1c_2...c_sW_s$ two decompositions
of the vector $(a_1,a_2,...,a_k)$. Then

\noindent i. \, $h=s$,

\noindent ii. The flag generated by the ordered free family
$(V_1,V_2,..,V_h)$ is equal to the flag generated by the ordered free
family $(W_1,W_2,...,W_h)$ that is $\forall i \in 1,..,h $
$$ \{
V_1,...,V_i\} =\{ W_1,...,W_i\} $$ where $\{U_i\}$ designates the linear space genrated by the vectors $U_i$.
\end{theorem}

\subsection{Geometrical interpretation of this decomposition}
Let $A$ be an $\mathbb{R}$ algebra of valuation. Consider a differential curve $\gamma$ in $\mathbb{R}^3$. We can embed $\gamma$ in a differential curve
$$\Gamma : \mathbb{R} \otimes A \rightarrow \mathbb{R}^3 \otimes A.$$
Let $t=t_0 \otimes 1 + 1 \otimes \epsilon$ an parameter infinitely close to $t_0$, that is $\epsilon \in \frak{m}$. If $M$ corresponds to the point of $\Gamma$ of parameter 
$t$ and $M_0$ those of $t_0$, then the coordinates of the point $M-M_0$ in the affine space $\mathbb{R}^3 \otimes A$ are in $\mathbb{R} \otimes \frak{m}$.
In the flag  associated to the decomposition
of $M-M_0$ we can considere a direct orthonormal frame $(V_1,V_2,V_3)$. It is the Serret-Frenet frame to $\gamma$ at the point $M_0$.

\section{Decomposition of a valued deformation of a Lie algebra}
\subsection{Valued deformation of Lie algebras}
Let $\frak{g}'_A$ be a valued deformation with base $A$ of the $\Bbb{K}$-Lie algebra $\frak{g}$. By definition, for every $X$ and $Y$ in $\frak{g}$
we have $[X,Y]_{\frak{g}'_A} -[X,Y]_{\frak{g}^{\,}_A} \in \frak{g} \otimes \frak{m}$. Suppose that $\frak{g}$ is finite dimensional and let $\{X_1,...,X_n \}$
be a basis of $\frak{g}$. In this case
$$[X_i,X_j]_{\frak{g}'_A} -[X_i,X_j]_{\frak{g}^{\,}_A}=\sum_{k} C_{ij}^k X_k$$
with $C_{ij}^k \in \frak{m}$. Using the decomposition of the vector of $\frak{m}^{n^2(n-1)/2}$ with for components $C_{ij}^k$ , we deduce that
$$
\begin{array}{lll}
\lbrack X_i,X_j]_{\frak{g}'_A} -[X_i,X_j]_{\frak{g}^{\,}_A}&=&a_{ij}(1)\phi _1(X_i,X_j)+a_{ij}(1)a_{ij}(2)\phi _2(X_i,X_j) \\
& &+...+
a_{ij}(1)a_{ij}(2)...a_{ij}(l)\phi _l(X_i,X_j)
\end{array}
$$
where $a_{ij}(s) \in \frak{m}$ and $\phi _1,...,\phi _l$ are linearly independent. The index $l$ depends of $i$ and $j$. Let $k$ be the supremum of indices $l$
when $1 \leq i, \ j \leq n$. Then we have 
$$
\begin{array}{lll}
\lbrack X,Y]_{\frak{g}'_A} -[X,Y]_{\frak{g}^{\,}_A}&=&\epsilon _1(X,Y)\phi _1(X,Y)+\epsilon _1(X,Y)\epsilon _2(X,Y)\phi _2(X,Y) \\
& &+...+
\epsilon _1(X,Y)\epsilon _2(X,Y)...\epsilon _k(X,Y)\phi _k(X,Y)
\end{array}
$$
where the bilinear maps $\epsilon _i$ have values in $\frak{m}$ and linear maps $\phi _i : \frak{g} \otimes \frak{g} \rightarrow \frak{g}$
are linearly independent.

If $\frak{g}$ is infinite dimensional with a countable basis $\{X_n\} _{n \in \mathbb{N}}$ then the $\mathbb{K}$-vector space of
linear map $T_2 ^1= \{ \phi : \frak{g} \otimes \frak{g} \rightarrow \frak{g} \}$ also admits a countable basis. 
\begin{theorem}
If $\mu _{\frak{g}'_A}$ (resp.
$\mu _{\frak{g}^{\,}_A}$) is the  law of the Lie algebra  ${\frak{g}'_A}$ (resp. ${\frak{g}^{\,}_A}$) then
$$\mu _{\frak{g}_A'}-\mu _{\frak{g}^{\,}_A}=\sum _{i \in I} \epsilon _1 \epsilon _2 ...\epsilon _i \phi_i$$
where $I$ is a finite set of indices, $\epsilon _i : \frak{g} \otimes \frak{g} \rightarrow \frak{m}$ are linear maps and $\phi _i$'s are linearly
independent maps in $T_2^1$.
\end{theorem}
\subsection{Equations of valued deformations}
We will prove that the classical equations of deformation given by Gerstenhaber are still valid in the general frame of valued deformations. Neverless
we can prove that the infinite system described by Gerstenhaber and which gives the conditions to obtain a deformation, can be reduced to a system of finite rank.
Let
$$\mu _{\frak{g}_A'}-\mu _{\frak{g}^{\,}_A}=\sum _{i \in I} \epsilon _1 \epsilon _2 ...\epsilon _i \phi_i$$
be a valued deformation of $\mu$ (the bracket of $\frak{g}$). Then $\mu _{\frak{g}_A'}$ satisfies the Jacobi equations. Following Gerstenhaber we
consider the Chevalley-Eilenberg graded differential complex $\mathcal{C}(\frak{g},\frak{g})$ and the product $\circ$ defined by
$$(g_q \circ f_p)(X_1,...,X_{p+q})=\sum (-1)^{\epsilon (\sigma)} g_q(f_p(X_{\sigma (1)},...,X_{\sigma (p)}),X_{\sigma (p+1)},...,X_{\sigma (q)})$$
where $\sigma$ is a permutation of ${1,...,p+q}$ such that $\sigma (1) < ...<\sigma (p)$ and $\sigma (p+1)<...<\sigma (p+q)$ (it is a 
$(p,q)$-schuffle);  $g_q \in
\mathcal{C}^q(\frak{g},\frak{g})$ and $f_p \in \mathcal{C}^p(\frak{g},\frak{g})$. As $\mu _{\frak{g}_A'}$ satisfies the Jacobi indentities,
$\mu _{\frak{g}_A'} \circ \mu _{\frak{g}_A'} = 0$. This gives
$$ (\mu _{\frak{g}_A}+\sum _{i \in I} \epsilon _1 \epsilon _2 ...\epsilon _i \phi_i) \circ ( \mu _{\frak{g}_A}+
\sum _{i \in I} \epsilon _1 \epsilon _2 ...\epsilon _i \phi_i)= 0. \quad \quad (1)$$
As $\mu _{\frak{g}_A} \circ \mu _{\frak{g}_A} = 0$, this equation becomes :
$$\epsilon _1(\mu _{\frak{g}_A} \circ \phi _1 + \phi _1 \circ \mu _{\frak{g}_A})+ \epsilon _1 U=0$$
where $U$ is in $\mathcal{C}^3(\frak{g},\frak{g}) \otimes \frak{m}$. If we symplify by $\epsilon_1$ which is supposed non zero if not the deformation is trivial, we obtain
$$(\mu _{\frak{g}_A} \circ \phi _1 + \phi _1 \circ \mu _{\frak{g}_A})(X,Y,Z)+  U(X,Y,Z)=0 $$
for all $X,Y,Z \in \frak{g}$. As $U(X,Y,Z)$ is in the module $\frak{g} \otimes \frak{m}$ and the first part in $\frak{g} \otimes A$, each one of these
vectors is null. Then $$(\mu _{\frak{g}_A} \circ \phi _1 + \phi _1 \circ \mu _{\frak{g}_A})(X,Y,Z)=0.$$
\begin{proposition}
For every valued deformation with base $A$ of the $\mathbb{K}$-Lie algebra $\frak{g}$, the first term $\phi$ appearing in the associated decomposition is a 2-cochain of the Chevalley-Eilenberg cohomology
of $\frak{g}$ belonging to $Z^2(\frak{g},\frak{g})$.
\end{proposition}
We thus rediscover the classical result of Gerstenhaber but in the broader context of valued deformations and not only for the valued deformation
of basis the ring of formal series.

In order to describe the properties of other terms of equations (1) we use the super-bracket of Gerstenhaber which endows the space of
Chevalley-Eilenberg cochains $\mathcal{C}(\frak{g},\frak{g})$ with a Lie superalgebra structure. When $\phi _i \in \mathcal{C}^2(\frak{g},\frak{g}),$
it is defines by
$$[\phi _i,\phi _j]=\phi _i \circ \phi _j + \phi _j \circ \phi _i $$
and $[\phi _i,\phi _j] \in \mathcal{C}^3(\frak{g},\frak{g})$.
\begin{lemma}
Let us suppose that $I=\{1,...,k\}$. If
$$\mu _{\frak{g}_A'}=\mu _{\frak{g}_A}+\sum _{i \in I} \epsilon _1 \epsilon _2 ...\epsilon _i \phi_i$$
is a valued deformation of $\mu$, then the 3-cochains $[\phi _i,\phi _j]$ and $[\mu, \phi_i]$, $1 \leq i,j \leq k-1$, generate 
a linear subspace $V$ of  $\mathcal{C}^3(\frak{g},\frak{g})$
of dimension less or equal to $k(k-1)/2$. Moreover, the 3-cochains $[\phi _i,\phi _j]$, $1 \leq i,j \leq k-1$, form a system of generators of
this space.
\end{lemma}

\noindent{\it Proof.}  Let $V$ be the subpace of $\mathcal{C}^3(\frak{g},\frak{g})$ generated by $[\phi _i,\phi _j]$ and $[\mu, \phi_i]$. If $\omega$ is
 a linear form on $V$ of which kernel contains the vectors $[\phi _i,\phi _j]$ for $1 \leq i,j \leq (k-1)$, then the equation (1) gives
$$\epsilon _1 \epsilon _2...\epsilon _k \omega ([\phi _1,\phi _k])+\epsilon _1 \epsilon _2^2...\epsilon _k \omega ([\phi _2,\phi _k])+...+
\epsilon _1 \epsilon _2^2...\epsilon _k^2 \omega ([\phi _k,\phi _k])+\epsilon _2\omega ([\mu,\phi _2]) $$
$$ +\epsilon _2\epsilon _3 \omega ([\mu,\phi _3])...+\epsilon _2\epsilon _3...\epsilon _k \omega ([\mu,\phi _k])=0.$$
As the coefficients which appear in this equation are each one in one $\frak{m}^p$, we have necessarily
$$\omega ([\phi _1,\phi _k])=...=\omega ([\phi _k,\phi _k])=\omega ([\mu,\phi _2])=...=\omega ([\mu,\phi _k])=0$$
and this for every linear form $\omega$ of which kernel contains $V$. This proves the lemma.

From this lemma and using the descending sequence
$$\frak{m} \supset \frak{m}^{(2)} \supset ... \supset \frak{m}^{(p)} ...$$
where $\frak{m}^{(p)}$ is the ideal generated by the products ${ a_1a_2...a_p, \  a_i \in \frak{m} }$ of length $p$,  we obtain :
\begin{proposition}
If
$$\mu _{\frak{g}_A'}=\mu _{\frak{g}_A}+\sum _{i \in I} \epsilon _1 \epsilon _2 ...\epsilon _i \phi_i$$
is a valued deformation of $\mu$, then we have the following linear system :
$$
\left\{
\begin{array}{l}

\delta \phi _2 = a_{11}^2 [\phi _1,\phi _1] \\
\delta \phi _3 = a_{12}^3 [\phi _1,\phi _2]+a_{22}^3[\phi _1,\phi _1]\\
... \\
\delta \phi _k = \sum _{1 \leq i \leq j \leq k-1} a_{ij}^{k} [\phi _i,\phi _j] \\
\lbrack \phi _1,\phi _k]=\sum _{1 \leq i \leq j \leq k-1} b_{ij}^{1} [\phi _i,\phi _j] \\
.... \\
\lbrack \phi _{k-1},\phi _k]=\sum _{1 \leq i \leq j \leq k-1} b_{ij}^{k-1} [\phi _i,\phi _j]

\end{array}
\right.
$$
where $\delta \phi _i = [\mu,\phi _i]$ is the coboundary operator of the Chevalley cohomology of the Lie algebra $\frak{g}$.
\end{proposition}
Let us suppose that the dimension of $V$ is the maximum $k(k-1)/2$. In this case we have no other relations between the generators
of $V$ and the previous linear system is complete, that is the equation of deformations does not give other relations than the relations of this system. 
The following result shows that, in this case, such deformation is isomorphic ,as Lie algebra laws,to a "polynomial"
valued deformation.
\begin{proposition}
Let be $\mu _{\frak{g}_A'}$ a valued deformation of $\mu$ such that
$$\mu _{\frak{g}_A'}=\mu _{\frak{g}_A}+\sum _{i=1,...,k} \epsilon _1 \epsilon _2 ...\epsilon _i \phi_i$$
and dim$V$=k(k-1)/2. Then there exists an automorphism of $\mathbb{K}^n \otimes \frak{m}$ of the form $f=Id \otimes P_k (\epsilon)$ with
$P_k(X) \in \mathbb{K}^{k}[X]$ satisfying $P_k(0)=1$ and $\epsilon \in \frak{m}$ such that the valued deformation $\mu _{\frak{g}''_A}$ defined by
$$\mu _{\frak{g}''_A}(X,Y)=h^{-1}(\mu _{\frak{g}'_A}(h(X),h(Y)))$$
is of the form
$$\mu _{\frak{g}_A"}=\mu _{\frak{g}_A}+\sum _{i=1,...,k} \epsilon ^i \varphi_i$$
where $\varphi _i = \sum _{j \leq i} \phi _j.$
\end{proposition}
\noindent {\it Proof.
 } Considering the Jacobi equation
$$ \lbrack\mu _{\frak{g}_A'},\mu _{\frak{g}_A'}]=0$$
and writting that dim$V$=$k(k-1)/2$, we deduce that there exist polynomials $P_i(X) \in \mathbb{K}[X]$ of degree $i$ such that
$$\epsilon _i = a_i \epsilon _k \frac{P_{k-i}(\epsilon _k)}{P_{k-i+1}(\epsilon _k)}$$ with $a_i \in \mathbb{K}$.
Then we have
$$\mu _{\frak{g}_A'}=\mu _{\frak{g}_A}+\sum _{i=1,...,k} a_1a_2...a_i (\epsilon _k)^i \frac{P_{k-i}(\epsilon _k)}{P_{k}(\epsilon _k)} \phi_i.$$
Thus
$$ P_{k}(\epsilon _k)\mu _{\frak{g}_A'}=P_{k}(\epsilon _k)\mu _{\frak{g}_A}+\sum _{i=1,...,k} a_1a_2...a_i (\epsilon _k)^i P_{k-i}(\epsilon _k) \phi_i.$$
If we write this expression according the increasing powers we obtain the announced expression. $\quad \Box$

\noindent Let us note that, for such deformation we have
$$
\left\{
\begin{array}{l}

\delta \varphi _2 + [\varphi _1,\varphi _1] =0\\
\delta \varphi _3 + [\varphi _1,\varphi _2]=0\\
... \\
\delta \varphi _k + \sum _{i+j = k}  [\varphi _i,\varphi _j] =0\\
\sum _{i+j=k+s}  [\varphi _i,\varphi _j] =0.\\

\end{array}
\right.
$$

\subsection{Particular case : one-parameter deformations of Lie algebras}
In this section the valuation ring $A$ is $\Bbb{K}[[t]]$. Its maximal ideal is $t\Bbb{K}[[t]]$ and the residual field is $\Bbb{K}$. 
Let $\frak{g}$ be
a $\mathbb{K}$- Lie algebra. Consider $ \frak{g} \otimes A$ as an $A$-algebra and let be $\frak{g}_A'$ a valued deformation of $\frak{g}$.
The bracket $[,]_t$ of this Lie algebra satisfies
$$[X,Y]_t=[X,Y] + \sum t^i \phi _i (X,Y).$$
Considered as a valued deformation wuth base $\mathbb{K}[[t]]$, this bracket can be written
$$[X,Y]_t = [X, Y] + \sum _{i=1}^{i=k}c_1(t)...c_i(t) \psi _i (X,Y)$$
where $(\psi _1,...,\psi _k)$ are linearly independent and $c_i(t) \in t\mathbb{C}[[t]]$. As $\phi _1 = \psi _1$, this bilinear map belongs
to $Z^2(\frak{g}, \frak{g})$ and we find again the classical result of Gerstenhaber. Let $V$ be the $\mathbb{K}$-vector space generated by
 $[\phi _i,\phi _j]$ and $[\mu, \phi _i], \ i,j=1,...,k-1$, $\mu$ being the law of $\frak{g}$. If dim$V =
k(k-1)/2$ we will say that one-parameter deformation $[,]_t$ is of maximal rank.
\begin{proposition}
Let
$$[X,Y]_t=[X,Y] + \sum t^i \phi _i (X,Y)$$
be a one-parameter deformation of $\frak{g}$. If its rank is maximal then this deformation is equivalent to a polynomial deformation
$$[X,Y]_t'=[X,Y]  + \sum _{i=1,...,k}t^i \varphi _i$$
with $\varphi _i =\sum _{j=1,...,i}a_{ij} \psi _j.$
\end{proposition}
\begin{corollary}
Every one-parameter deformation of maximal rank is equivalent to a local non valued deformation with base the local algebra $\mathbb{K}[t]$.
\end{corollary}
Recall  that the algebra $\mathbb{K}[t]$ is not an algebra of valuation. But every local ring is dominated by a valuation ring. Then this corollary can be interpreted as saying that every deformation in the local
algebra $\mathbb{C}[t]$ of polynomials with coefficients in $\mathbb{C}$ is equivalent to a "classical"-Gerstenhaber deformation with maximal rank.

\section{Deformations of the enveloping algebra of a rigid Lie algebra}
\subsection{Valued deformation of associative algebras}
Let us recall that the category of $\mathbb{K}$-associative algebras is a monoidal category. 
\begin{definition}
Let $\frak{a}$ be  a $\mathbb{K}$-associative algebra and $A$ an  $\mathbb{K}$-algebra of valuation of such that the residual field
$\frac{A}{\frak{m}}$ is isomorphic to $\mathbb{K}$ (or to a subfield $\Bbb{K}'$ of $\mathbb{K}$). A valued deformation of $\frak{a}$ with  
base $A$ is an $A$-associative algebra $\frak{a}'_A$
such that the underlying $A$-module of $\frak{a}'_A$ is $\frak{a}_A$ and that
$$(X.Y)_{\frak{a}'_A} -(X.Y)_{\frak{a}^{\,}_{A}} $$ belongs to the $\frak{m}$-quasi-module $\frak{a} \otimes \frak{m}$ where $\frak{m}$ is the maximal ideal of
$A$.
\end{definition}
The classical one-parameter deformation is a valued deformation. As in the Lie algebra case we can develop the decomposition of a valued
deformation. It is sufficient to change the Lie bracket by the associative product and the Chevalley cohomology by the Hochschild cohomology.

The most important example concerning valued deformations of associative algebras is those of the associative algebra of smooth fonctions of a manifold.
But we will be interested here by associative algebras that are the enveloping algebras of  Lie algebras. More precisely, what can we say about
the valued deformations of the enveloping algebra of a rigid Lie algebra?

\subsection{Complex rigid Lie algebras}
In this section we suppose that $\mathbb{K}=\mathbb{C}$.
Let ${\cal{L}}_n$ be the algebraic variety of structure constants of $n$-dimensional complex Lie algebra laws. The basis of ${\mathbb{C}}^n$ 
being fixed, we can identify a law with its structure constants. Let us consider the action of the linear group $Gl(n, \Bbb{C})$ on ${\cal{L}}_n$ :
$$ \mu'(X,Y) = f^{-1} \mu (f(X),f(Y)).$$
We denote by $\cal{O}({\mu}) $ the orbit of $\mu$.
\begin{definition}
The law $\mu \in {\cal{L}}_n$ is called rigid if $\cal{O}({\mu}) $ is Zariski-open in ${\cal{L}}_n$.
\end{definition}
Let $\frak{g}$ be a $n$-dimensional complex Lie algebra with product $\mu$ and $\frak{g}_A$ a valued deformation with base $A$. As before
$\mathcal{F}_A$ is the field of fractions of $A$.
\begin{definition}
Let $A$ be a valued $\Bbb{C}$-algebra. We say that $\frak{g}$ is $A$-rigid if for every valued deformation $\frak{g}_A'$ of $\frak{g}_A$ there exists
a $\mathcal{F}_A$-linear isomorphism between $\frak{g}_A'$ and $\frak{g}_A$.
\end{definition}

Let $\mu _{\frak{g}'_A}$ be a valued deformation of
$\mu_{\frak{g}^{\,}_A}$. If we write $\mu
_{\frak{g}'_A}-\mu_{\frak{g}^{\,}_A}=\phi$, then $\phi(X,Y) \in
\frak{g} \otimes \frak{m}$ for all $X,Y \in \frak{g} \otimes A.$ If
$\mu_{\frak{g}^{\,}_A}$ is rigid, there exits $f \in Gl_n(\frak{g} \otimes
\cal{F}_A)$ such that
$$f^{-1}(\mu_{\frak{g}'_A}(f(X),f(Y)))=\mu_{\frak{g}^{\,}_A}(X,Y).$$
Thus
$$\mu_{\frak{g}^{\,}_A}(f(X),f(Y))-f(\mu_{\frak{g}^{\,}_A}(X,Y))=\phi(f(X),f(Y)).$$
As $\frak{g}^{\,}_A$ is invariant by $f$, $\phi(f(X),f(Y)) \in
\frak{g} \otimes \frak{m}.$ So we can decompose $f$ as $f=f_1+f_2$
with $f_1 \in Aut(\frak{g}_A)$ and $f_2:\frak{g}_A \rightarrow
\frak{g} \otimes \frak{m}.$ Let $f'$ be $f'=f \circ f^{-1}_1.$ Then
$$f'^{-1}(\mu_{\frak{g}'_A}(f'(X),f'(Y)))=\mu_{\frak{g}^{\,}_A}(X,Y)$$
and $f'=Id+h$ with $h:\frak{g}_A \rightarrow\frak{g} \otimes
\frak{m}.$ Thus we have proved
\begin{lemma}
If $\mu_{\frak{g}^{\,}_A}$ is $A$-rigid for every valued deformation
$\mu_{\frak{g}'_A}$ there exits $f \in Gl_n(\frak{g} \otimes
\cal{F}_A)$ of the type $f=Id+h$ with $h:\frak{g}_A \rightarrow\frak{g} \otimes
\frak{m}$ such that
$$f^{-1}(\mu_{\frak{g}'_A}(f(X),f(Y)))=\mu_{\frak{g}^{\,}_A}(X,Y)$$
for every $X,Y \in \frak{g}_A.$
\end{lemma}

\bigskip

\noindent {\bf Remark.} If $f=Id+h$ then $f^{-1}=Id+k$. As
$\frak{g}_A$ is invariant by $f$, the linear map $k$ satisfies $k:
\frak{g}_A \rightarrow \frak{g} \otimes \frak{m}.$

\begin{theorem}
If the residual field of the valued ring is isomorphic to $\Bbb{C}$ then the notions of $A$-rigidity and of rigidity are equivalent.
\end{theorem}
{\it Proof.} Let us suppose that for every valued algebra of residual field $\mathbb{C}$, the Lie algebra $\frak{g}$ is $A$-rigid. We will consider
the following special valued algebra: let $\mathbb{C}^*$ be non standard extension of $\mathbb{C}$ in the Robinson sense ([Ro]). If $\mathbb{C}_l$
is the subring of non-infinitely large elements of  $\mathbb{C}^*$ then the subring $\frak{m}$ of infinitesimals is the maximal ideal of
$\mathbb{C}_l$ and $\mathbb{C}_l$ is a valued ring. Let us consider $A=\mathbb{C}_l$. In this case we have a natural embedding of the
variety of $A$-Lie algebras in the variety of $\mathbb{C}$-Lie algebras. Up this embedding (called the transfert principle in the Robinson theory), the set of $A$-deformations of $\frak{g}_A$ is an infinitesimal
neighbourhood of $\frak{g}$ contained in the orbit of $\frak{g}$. Then
$\frak{g}$ is rigid. $\quad \Box$
\medskip

\noindent {\bf Examples.} If
 $A=\Bbb{C}[[t]]$ then $\Bbb{K}'=\Bbb{C}$ and we find again the
classical approach to the rigidity. We have another example, yet
used in the proof of Theorem 3, considering
a non standard extension $\Bbb{C}^*$ of $\Bbb{C}$. In this context the notion of rigidity has been developed in [A.G] 
(such a deformation is called perturbation). This work
has allowed to classifiy complex finite dimensional rigid Lie algebras up the dimension eight.

\subsection{Deformation of the enveloping algebra of a Lie algebra}
Let $\frak{g}$ be a finite dimensional $\mathbb{K}$-Lie algebra and $\mathcal{U}(\frak{g})$ its enveloping algebra. In this section we consider
a particular valued deformation of $\U$ corresponding to the valued algebra $\mathbb{K}[[t]]$. In [P], the following result
is proved:
\begin{proposition}
If $\frak{g}$ is not rigid then $\U$ is not $\mathbb{K}[[t]]$-rigid.
\end{proposition}
Recall that if the Hochschild cohomology $H^*(\U,\U)$ of $\U$ satisfies $$H^2(\U,\U)=0$$ then $\U$ is $\mathbb{K}[[t]]$-rigid. By the Cartan-Eilenberg
theorem, we have that
$$H^n(\U,\U)=H^n(\frak{g},\U).$$
\begin{theorem}[P]
Let $\frak{g}$ be a  rigid Lie algebra. If $H^2(\frak{g},\mathbb{C}) \neq 0$, then $\U$ is not $\mathbb{K}[[t]]$-rigid.
\end{theorem}
From [C] and [A.G] every solvable complex Lie algebra decomposes as $\frak{g}=\frak{t} \oplus \frak{n}$ where $\frak{n}$ is the niladical of
$\frak{g}$ and $\frak{t}$ a maximal  exterior torus of derivations in the Malcev sense. Recall that the rank of $\frak{g}$ is the dimension of
$\frak{t}$. A direct consequence of Petit's theorem is that for every complex rigid Lie algebra of rank equal or greater than 2 its envelopping algebra
is not rigid.
\begin{theorem}
Let $\frak{g}$ be a complex finite dimensional rigid Lie algebra of rank $1$. Then
$$dim \ H^2(\frak{g},\mathbb{C}) =0$$
if and only if $0$ is not a root of the nilradical $\frak{n}.$
\end{theorem}
{\it Proof}. Suppose first that $0$ is not a root of
$\frak{n}$ that is for every $X \neq 0  \in  \frak{t}$, 0 is not an eigenvalue of the semisimple
operator $ad X$. Let $\theta$ be in $Z^2(\frak{g},\mathbb{C})$. Let $(X,Y_i)_{i=1,...,n-1}$ a basis of $\frak{n}$ adapted to the decomposition
$\frak{g}=\frak{t} \oplus \frak{n}$. In particular we have
$$[X,Y_i]=\lambda _i Y_i$$
with $\lambda _i \in \mathbb{N}^*$ for all $i=1,..,n-1$ ([A.G]). As
$d\theta =0$ we have for all $i,j=1,...,n-1$
$$d\theta(X,Y_i,Y_j)=\theta(X,[Y_i,Y_j])+\theta(Y_i,[Y_j,X])+\theta(Y_j,[X,Y_i])=0,$$
for all $1 \leq i,j \leq k-1$, and this gives
$$(\lambda_i+\lambda_j)\theta(Y_i,Y_j)= \theta(X,[Y_i,Y_j]). \quad
\quad \quad (1)$$
If $(\lambda_i+\lambda_j)$ is not a root, then $[Y_i,Y_j]=0$ and this implies that $\theta(Y_i,Y_j)= 0$. If not, $(\lambda_i+\lambda_j)=\lambda _k$
is a root. Let us note ${Y_k^1,...,Y_k^{n_k}}$ the eigenvectors of the chosen basis corresponding to the root $\lambda_k$. We have
 $$[Y_i,Y_j]=\sum _{s=1}^{n_k} a_{ij}^s(k) Y_k^s.$$
Let us consider the dual basis $\{\omega _0,\omega _1,...,\omega _{n-1}\}$ of $\{X,Y_1,...,Y_{n-1}\}$. We have
$$d\omega _k^s = \lambda _k \omega _0 \wedge \omega _k^s + \sum _{l,m} a_{lm}^s(k) \omega _l \wedge \omega _m $$
where the pairs $(l,m)$ are such that $\lambda _l + \lambda _m = \lambda _k.$
Then we deduce from (1)
$$\sum a_{ij}^s(k) \theta (X,Y_k^s) - \lambda _k\theta(Y_i,Y_j)=0.$$
Let us fix $\lambda_k$. If we write
$$
\begin{array}{ll}
\theta = & \sum _{l,m, \lambda _l + \lambda _m =\lambda _k} b_{lm}(k)\w _l \wedge \w _m +
\sum _{r,s, \lambda _r + \lambda _s \neq \lambda _k}c _{rs}(k)
 \w _r \wedge \w _s  \\
& \\
 &
 + \sum _k \sum _{s=1}^{n_k} \beta _k^s \w _0 \wedge \w _k^s
\end{array}
$$
then, for every pair $(i,j)$ such that
$\lambda_i+\lambda_j=\lambda_k$, (1) gives
$$-\lambda _k b_{ij}(k) + \sum _{s=1}^{n_k}  a_{ij}^s(k) \beta _{k}^s = 0$$
and
$$b_{ij}(k)= \sum _{s=1}^{n_k} \frac{a_{ij}^s(k)}{\lambda _k}\beta _{k}^s.$$
The expression of $\theta$ becomes
$$
\begin{array}{ll}
\theta = & \sum _k \sum _{s=1}^{n_k} \beta _k^s \w _0 \wedge \w _k^s +
\sum _{i,j, \lambda _i + \lambda _j =\lambda _k} \frac{a_{ij}^s(k)}{\lambda _k}\beta _{k}^s \w _i \wedge \w _j  \\
 & \\
& + \sum _{r,s, \lambda _r + \lambda _s \neq \lambda _k}c _{rs}(k)
 \w _r \wedge \w _s .
\end{array}
$$
Thus
$$
\begin{array}{lll}
\theta & = & \sum _k (\frac{1}{\lambda _k} \sum _{s=1}^{n_k} \beta _k^s ( \lambda _k \w _0 \wedge \w _k^s +
\sum _{i,j, \lambda _i + \lambda _j =\lambda _k} a_{ij}^s(k) \w _i \wedge \w _j ) \\
 & & \\
& & + \sum _{r,s, \lambda _r + \lambda _s \neq \lambda _k}c _{rs}(k)
 \w _r \wedge \w _s \\
& \\
& = & \sum _s \beta _k^s d\w _k^s + \sum _{k' \neq k} \sum _{s=1}^{n_{k'}} \beta _{k'}^s \w _0 \wedge \w _{k'}^s \\
& & \\
& & +
\sum _{r,s, \lambda _r + \lambda _s \neq \lambda _k}c _{rs}(k)\w _r \wedge \w _s .
\end{array}
$$
If we continue this method for all the non simple roots (that is which admit a decomposition as sum of two roots, we obtain
the heralded result.

\noindent For the converse, if $0$ is a root, then the cocycle
$$ \theta=\omega_0 \wedge \omega '_0$$
where $\omega '_0$ is related with the eigenvector associated to the
root $0$ is not integrable.  $\quad \Box$
\medskip

\noindent {\bf Remark.}
It is easy to verify that every solvable rigid Lie algebra of rank greater or equal to $2$ cannot have $0$ as root. Likewise every solvable rigid
Lie algebra of rank $1$ and of dimension less than $8$ has not $0$ as
root. This confirm in small dimension the following conjecture [Ca]:

\noindent {\it
If $\frak{g}$ is a complex solvable finite dimensional rigid Lie algebra of rank $1$, then $0$ is not a root.
}
\medskip

\noindent{\bf Consequences.}
If $H^2(\frak{g},\Bbb{C})\neq 0$, there exits $\theta \in \wedge^2
\frak{g}^*$ such that $[\theta]_{H^2} \neq 0.$ If $rg(\frak{g})\geq 2$,
then we can suppose that $\theta \in \wedge^2\frak{t}^*$ and $\omega$
defines a non trivial deformation of $\mathcal{U} (\frak{g}).$
If $rg(\frak{g})=1$, then $0$ is not a root of $\frak{t}$. The Hochschild
Serre sequence gives:
$$
\begin{array}{lll}
H^2_{CE}(\frak{g},\mathcal{U}(\frak{g} ) ) & = & (\wedge ^2\frak{t}^*\otimes
Z( \mathcal{U}( \frak{g} ) ) ) \oplus ( \frak{t}^* \otimes H^1_{CE}(
\frak{n},\mathcal{U}( \frak{g} ) )^t )\oplus H^2_{CE}(
\frak{n},\mathcal{U}(\frak{g} ) )^t \\
& = & \frak{t}^* \otimes H^1_{CE}(\frak{n},\mathcal{U}( \frak{g}))^t
\oplus H^2_{CE}( \frak{n},\mathcal{U}(\frak{g}) )^t
\end{array}
$$
But from the previous proof, if $\theta$ is a non trivial 2-cocycle of
$Z^2_{CE}(\frak{g},\Bbb{C} )$ then $i(X) \theta \neq 0$ for every $X \in
\frak{t}$, $X \neq 0.$ The 1-form $\omega = i(X) \theta$ is closed. Then
$\theta$ corresponds to a cocycle belonging to $\frak{t}^* \otimes
Z^1_{CE}( \frak{n},\mathcal{U}(\frak{g}))^t$ and defines a deformation
of $\mathcal{U}(\frak{g}). $

\begin{theorem}
Let $\frak{g}$ be a solvable complex rigid Lie algebra. If its rank is greater or equal to $2$
or if the rank is $1$ and $0$ is a root, then the enveloping algebra $\cal{U}(\frak{g})$ is not rigid.
\end{theorem}

\noindent {\bf Remark.} In [P], T.Petit discribes some examples of deformations of
the enveloping algebra of a rigid Lie algebra $\frak{g}$ in small dimension and satisfying
$H^2_{CE}(\frak{g},\Bbb{C}) = 0.$ For this, he shows that every
deformation of the linear Poisson structure on the dual $\frak{g}^*$
of $\frak{g}$ induces a non trivial deformation of
$\mathcal{U}(\frak{g} ).$ This reduces the problem to find non trivial
deformation of the linear Poisson structure.

\subsection{Poisson algebras}
Recall that a Poisson algebra $\cal{P}$ is a (commutatitive) associative algebra endowed with a second algebra law satisfying
the Jacobi's identity and the Liebniz rule
$$[a,bc]=b[a,c]+[a,b]c$$
for all $a,b,c \in \cal{P}$.
The tensor product ${\cal{P}}_1 \otimes {\cal{P}}_2$ of two Poisson algebras is again a Poisson algebra with the following associative
and Lie products on ${\cal{P}}_1 \otimes {\cal{P}}_2 :$
$$(a_1 \otimes a_2).(b_1 \otimes b_2)=(a_1.b_1) \otimes (a_2.b_2)$$
$$[(a_1 \otimes a_2),(b_1 \otimes b_2) ]=([a_1,b_1] \otimes a_2.b_2 )+(a_1.b_1 \otimes [a_2,b_2])$$
for all $a_1,b_1 \in {\cal{P}}_1, a_2,b_2 \in {\cal{P}}_2.$ We can verify easily that these laws satisfy the Leibniz rule.

Every commutative associative algebra has a natural Poisson structure, putting $[a,b]=ab-ba=0$. Then the tensor product of a Poisson algebra
by a valued algebra is as well a Poisson algebra. In this context we have the notion of valued deformation. For example, if we take as valued algebra
the algebra $\Bbb{C}[[t]]$, then the Poisson structure of ${\cal{P}} \otimes {\Bbb{C}}[[t]]$ is given by
$$(a_1 \otimes a_2(t)).(b_1 \otimes b_2(t))=(a_1.b_1) \otimes (a_2(t).b_2(t))$$
$$[(a_1 \otimes a_2(t)),(b_1 \otimes b_2(t)) ]=[a_1,b_1] \otimes a_2(t).b_2(t) $$
because $\Bbb{C}[[t]]$ is a commutative associative algebra.

\medskip

\noindent{\bf Remark.} As we have a tensorial category it is natural to look if we can define a Brauer Group for Poisson algebras. As the associative
product corresponds to the classical tensorial product of associative algebras, we can consider only Poisson algebras which are finite dimensonal
simple central algebras. The matrix algebras $M_n(\Bbb{C})$ are Poisson algebras. Then, considering the classical equivalence relation for define the Brauer
Group, the class of matrix algebra constitutes an unity. Now the opposite algebra $A^{op}$ also is a Poisson algebra. In fact the associative product
is given by $a._{op}b=ba$ and the Lie bracket by $[a,b]_{op}=ba-ab$. Thus
$$[a,b._{op}c]_{op}=[a,cb]_{op}=cba-acb$$
$$b._{op}[a,c]_{op}+[a,b]_{op}._{op}c=cab-acb+cba-cab=cba-acb$$
this gives the Poisson structure of $A^{op}$. The opposite algebra $A^{op}$ is, modulo the equivalence relation, the inverse of $A$.

\section{Deformations of non associative algebras}
\subsection{Lie-admissible algebras}
In [R], special classes of non-associative algebras whose
laws give a Lie bracket by anticommutation are presented. If $A$ is a
$\Bbb{K}$-algebra,  we'll note by $a_{\mu}$ the associator of its law
$\mu$:
$$a_{\mu}(X,Y,Z)=\mu(\mu(X,Y),Z)-\mu(X,\mu(Y,Z)).$$
Let $\Sigma_n$ be the $n$-symmetric group.
\begin{definition}
An algebra $A$ is Lie-admissible if
$$ \sum_{\sigma \in \Sigma_3}a_{\mu} \circ \sigma =0,$$
with $\sigma(X_1,X_2,X_3)=(X_{\sigma^{-1}(1)},X_{\sigma^{-1}(2)},X_{\sigma^{-1}(3)}).$
\end{definition}
Let $G$ be a sub-group of $\Sigma_3$. The Lie-admissible algebra
$(A,\mu)$ is called $G$-associative if
$$ \sum_{\sigma \in G}a_{\mu} \circ \sigma =0.$$
Let us note that this last identity implies the Lie-admissible
identity. If $G$ is the trivial sub-group, then the
corresponding class of $G$-associative algebras is nothing other that
the associative algebras but for all other sub-group, we have
non-assotiative algebras. For example, if $G=<Id,\tau_{23} >$ we
obtain the pre-Lie algebras ([G]). If $G=<Id,\tau_{12}>$ the
corresponding algebras are the Vinberg algebras.
\subsection{External tensor product}
It is easy to see that each one of the categorie of $G$-associative
algebras is not tensorial exept for $G=<Id>.$ But in [G.R] we have
proved the following result:
\begin{theorem}
For every sub-group $G$ of $\Sigma_3$, let $\cal{G}-\cal{A}$$ss$ be the
associated operad and $\cal{G}-\cal{A}$$ss^{!}$ its dual operad. For
every $\cal{G}-\cal{A}$$ss$-algebra $A$ and
$\cal{G}-\cal{A}$$ss^{!}$-algebra $B$, $A \otimes B$ is a $
\cal{G}-\cal{A}$$ss$-algebra.
\end{theorem}
Recall the structure of $\cal{G}-\cal{A}$$ss^{!}$-algebras.

-If $G=<Id>$, then the $\cal{G}-\cal{A}$$ss^{!}$ are the associative
 algebras,

-If $G=<Id,\tau_{12}>$, then the $\cal{G}-\cal{A}$$ ss^{!}$ are the associative
 algebras satisfying $abc=bac$,

-If $G=<Id,\tau_{23}>$, then the $\cal{G}-\cal{A}$$ ss^{!}$ are the associative
 algebras satisfying $abc=acb$,

-If $G=<Id,\tau_{13}>$, then the $\cal{G}-\cal{A}$$ ss^{!}$ are the associative
 algebras satisfying $abc=cab$,

-If $G=\cal{A}_3$, then the $\cal{G}-\cal{A}$$ ss^{!}$ are the associative
 algebras satisfying $abc=bac=cab$,

-If $G=\Sigma_3$, then the $\cal{G}-\cal{A}$$ ss^{!}$ are the
 3-commutative  associative
 algebras ([R]).
\subsection{Valued deformation of $G$-associative algebras}
\begin{definition}
Let $\cal{A}$ be a $\cal{G}-\cal{A}$$ss$ algebra and $A$ a valued
$\cal{G}-\cal{A}$$ss^{!}$ algebra. Let $\mu$ be the law of the
$A$-algebra $\cal{A} \otimes A$. A deformation of $\cal{A}$ of basis
$A$ is a $A$-algebra whose law $\mu'$ satisfies:
$$\mu'(X,Y)-\mu(X,Y) \in \cal{A} \otimes \frak{m}$$
for every $X,Y \in \cal{A}$, where $\frak{m}$ is the maximal two-sided
ideal of $A$.
\end{definition}
If $A$ is a commutative associative algebra, then it's a
$\cal{G}-\cal{A}$$ss^{!}$-algebra for every $G$. We can study the
valued deformations in this particular case.

\bigskip

\noindent {\bf Acknowledgements} We express our thanks to Martin Markl
for reading the manuscript and many helpful remarks.

\bigskip

\noindent {\bf References}

\medskip

\noindent [A.G] Ancochea Bermudez J.M, Goze M., On the classification of Rigid Lie algebras.{\it Journal of algebra}, {\bf 245} (2001), 68-91.

\medskip

\noindent [C] Carles R., Sur la structure des alg\`ebres de Lie rigides. {\it Ann. Inst. Fourier}, {\bf 34} (1984), 65-82.

\medskip

\noindent [Ca] Campoamor R., Invariants of solvable rigid Lie algebras up to dimension 8. {\it J. Phys. A: Math Gen.}, {\bf 35} (2002), 6293-6306.

\medskip

\noindent [F] Fialowski A., Post, G., Versal deformation of the Lie algebra $L\sb 2$. {\it J. Algebra}  {\bf 236},1,  (2001), 93-109.

\medskip

\noindent [G] Gerstenhaber M., On the deformations of ring and algebras, {\it Ann. Math.} {\bf 74}, 11, (1964), 59-103.

\medskip

\noindent [G.R] Goze M., Remm E., On the categories of Lie-admissible. Preprint Mulhouse (2001) xxx RA/0210291
algebras,

\medskip

\noindent [L.S] Lichtenbaum S., Schlessinger M., The cotangent complex of a morphism.  Trans. Amer. Math. Soc.  128  1967 41-70.

\medskip

\noindent [M.S] Markl, Martin; Stasheff, James D. Deformation theory via deviations. {\it J. Algebra} 
 {\bf 170}, 1,  (1994),  122-155. 

\medskip
 
\noindent [P] Petit T., Alg\`ebres enveloppantes des alg\`ebres de Lie
rigides. Th\`ese Mulhouse 2001.

\medskip

\noindent [R] Remm E., Op\'erades Lie-admissibles. {\it
C.R.Acad. Sci. Paris}, ser I 334 (2002) 1047-1050.

\medskip

\noindent [S] Schlessinger M., Functors of Artin rings.  Trans. Amer. Math. Soc.  130  1968 208-222.

\end{document}